\documentclass[12pt]{article}
\usepackage{latexsym,amssymb,amsmath}
\usepackage{stmaryrd}
\begin{document}
\baselineskip=20pt

\newcommand{\la}{\langle}
\newcommand{\ra}{\rangle}
\newcommand{\psp}{\vspace{0.4cm}}
\newcommand{\pse}{\vspace{0.2cm}}
\newcommand{\ptl}{\partial}
\newcommand{\dlt}{\delta}
\newcommand{\sgm}{\sigma}
\newcommand{\al}{\alpha}
\newcommand{\be}{\beta}
\newcommand{\G}{\Gamma}
\newcommand{\gm}{\gamma}
\newcommand{\vs}{\varsigma}
\newcommand{\Lmd}{\Lambda}
\newcommand{\lmd}{\lambda}
\newcommand{\td}{\tilde}
\newcommand{\vf}{\varphi}
\newcommand{\yt}{Y^{\nu}}
\newcommand{\wt}{\mbox{wt}\:}
\newcommand{\rd}{\mbox{Res}}
\newcommand{\ad}{\mbox{ad}}
\newcommand{\stl}{\stackrel}
\newcommand{\ol}{\overline}
\newcommand{\ul}{\underline}
\newcommand{\es}{\epsilon}
\newcommand{\dmd}{\diamond}
\newcommand{\clt}{\clubsuit}
\newcommand{\vt}{\vartheta}
\newcommand{\ves}{\varepsilon}
\newcommand{\dg}{\dagger}
\newcommand{\tr}{\mbox{Tr}}
\newcommand{\ga}{{\cal G}({\cal A})}
\newcommand{\hga}{\hat{\cal G}({\cal A})}
\newcommand{\Edo}{\mbox{End}\:}
\newcommand{\for}{\mbox{for}}
\newcommand{\kn}{\mbox{ker}}
\newcommand{\Dlt}{\Delta}
\newcommand{\rad}{\mbox{Rad}}
\newcommand{\rta}{\rightarrow}
\newcommand{\mbb}{\mathbb}
\newcommand{\lra}{\Longrightarrow}

\title{\bf  Generalized Conformal Representations
of Orthogonal Lie Algebras  \footnote{Zhao's research is supported
by NSFC Grant 10701002; Xu's research is supported by NSFC Grant
10871193.}}

\author{$\mbox{Xiaoping Xu}^1$  and $\mbox{Yufeng Zhao}^2$} \maketitle

{\small \noindent 1. Corresponding author, Hua Loo-Keng Key
Mathematical Laboratory, Institute of Mathematics, Academy of
Mathematics and Systems Sciences, Chinese
Academy of  Sciences, Beijing 100190, P. R. China.\\
2. LMAM, School of Mathematical Sciences, Peking University,
  Beijing 100871, P. R. China}

\begin{abstract}

\quad The conformal transformations with respect to the metric
defining $o(n,\mbb{C})$ give rise to a nonhomogeneous polynomial
representation of $o(n+2,\mbb{C})$. Using Shen's technique of mixed
product, we generalize the above representation to a non-homogenous
representation of $o(n+2,\mbb{C})$ on the tensor space of any
finite-dimensional irreducible $o(n,\mbb{C})$-module with the
polynomial space, where a hidden central transformation is involved.
Moreover, we find a condition on the constant value taken by the
central transformation  such that the generalized conformal
representation is irreducible. In our approach,  Pieri's formulas,
invariant operators and the idea of Kostant's characteristic
identities play key roles. The result could be useful in
understanding higher-dimensional
 conformal field theory with the constant value taken by the central transformation as the central charge.
Our representations virtually provide natural
  extensions of the conformal transformations on
 a Riemannian manifold to its vector bundles.

  \end{abstract}

\section {Introduction}

 \quad \quad A quantum field is an operator-valued function on a certain Hilbert space,
which is often a direct sum of  infinite-dimensional irreducible
modules of a certain Lie algebra (group). The Lie algebra of
two-dimensional conformal group is exactly the Virasoro algebra,
which is infinite-dimensional. The minimal models of two-dimensional
conformal field theory were constructed from direct sums of certain
infinite-dimensional irreducible modules of the Virasoro algebra,
where a distinguished module called, the {\it vacuum module}, gives
rise to a vertex operator algebra. When $n\geq 3$,  $n$-dimension
conformal groups (depending on the metric) are finite-dimensional.
Higher-dimensional conformal field theory is not so well understood
partly because we are lack of enough knowledge on the
infinite-dimensional irreducible modules of orthogonal Lie algebras
 that are compatible to the natural
conformal representations. This motivates us to study explicit
infinite-dimensional irreducible representations of the orthogonal
Lie algebra $o(n+2,\mbb{C})$ by using the non-homogeneous polynomial
representations arising from conformal transformations with respect
to the metric defining $o(n,\mbb{C})$ and Shen's technique of mixed
product (cf. [Sg1-3]) (also known as Larsson functor (cf. [La])).

It is well known that $n$-dimensional projective group gives rise to
a non-homogenous representation of the Lie algebra $sl(n+1,\mbb{C})$
on the polynomial functions of the projective space. Using Shen's
mixed product for Witt algebras, the  authors [ZX] generalized the
above representation of $sl(n+1,\mbb{C})$ to a non-homogenous
representation on the tensor space of any finite-dimensional
irreducible $gl(n,\mbb{C})$-module with the polynomial space.
Moreover, the structure of such a representation was completely
determined by employing projection operator techniques (cf. [Gm1])
and the well-known Kostant's characteristic identities (cf. [K]).

 In this paper, we generalize the conformal representation of of $o(n+2,\mbb{C})$ to a non-homogenous
representation of $o(n+2,\mbb{C})$ on the tensor space of any
finite-dimensional irreducible $o(n,\mbb{C})$-module with the
polynomial space by Shen's idea of mixed product for Witt algebras.
It turns out that  a hidden central transformation is involved. More
importantly,  we find a condition on the constant value taken by the
central transformation  such that the generalized conformal
representation is irreducible. In our approach,  Pieri's formulas,
invariant operators and the idea of Kostant's characteristic
identities play key roles. The result could be useful in
understanding higher-dimensional
 conformal field theory with the constant value taken by
 the central transformation as the central charge. Our representations virtually provide natural
  extensions of the conformal transformations on
 a Riemannian manifold to its vector bundles.

Characteristic identities have a long history. The first person to
exploit them was Dirac [D], who wrote down what amounts to the
characteristic identity for the Lie algebra $o(1,3)$. This
particular example is intimately connected with the problem of
describing the structure of  relativistically invariant wave
equations. It had been shown by Kostant [K] (also cf. [Gm2]) that
the characteristic identities for semi-simple Lie algebras also hold
for infinite dimensional representations.

The $n$-dimensional conformal group with respect to Euclidean metric
$(\cdot,\cdot)$ is generated by the translations, rotations,
dilations and special conformal transformations
$$\vec x\mapsto\frac{\vec x-(\vec x,\vec x)\vec b}{(\vec b,\vec b)
(\vec x,\vec x)-2(\vec b,\vec x)+1}.\eqno(1.1)$$
   Let ${\cal
A}=\mbb{C}[x_1,...,x_n]$. The Witt algebra ${\cal
W}(n)=\sum_{i=1}^n{\cal A}\ptl_{x_i}$ with the commutator of
differential operators as its Lie bracket. The conformal
transformations with respect to the metric defining $o(n,\mbb{C})$
give rise to a nonhomogeneous polynomial representation of
$\vt:o(n+2,\mbb{C})\rta {\cal W}(n)$, acting on ${\cal A}$. Let
$E_{r,s}$ be the square matrix with 1 as its $(r,s)$-entry and 0 as
the others. Acting on the entries of the elements of the Lie algebra
$gl(n,{\cal A})$, ${\cal W}(n)$ becomes a Lie subalgebra of the
derivation Lie algebra  of $gl(n,{\cal A})$. In particular,
$$\widehat{\cal W}(n)={\cal W}(n)\oplus gl(n,{\cal A})\eqno(1.2)$$
becomes a Lie algebra with the Lie bracket
$$[d_1+A_1,d_2+A_2]=[d_1,d_2]+[A_1,A_2]+d_1(A_2)-d_2(A_1)\eqno(1.3)$$
for $d_1,d_2\in{\cal W}(n)$ and $A_1,A_2\in gl(n,{\cal A})$. Shen
[S1] found a monomorphism $\Im:{\cal W}(n)\rightarrow \widehat{\cal
W}(n)$ defined by
$$\Im(\sum\limits_{i=1}^{n}f_{i}\ptl_{x_i})=
\sum\limits_{i=1}^{n}f_{i}\ptl_{x_i}+\sum_{i,j=1}^n\ptl_{x_i}(f_j)\otimes
E_{i,j}.\eqno(1.4)$$ Moreover, he [S1-S3] used the monomorphism
$\Im$ to develop a theory of mixed products for the modules of Lie
algebras of Cartan type, which is also known as  the Larsson functor
(cf. [L]) in the case of Witt algebras. Rao [R] constructed some
irreducible weight modules over the derivation Lie algebra of the
algebra of Laurent polynomials based on Shen's mixed product. Lin
and Tan [LT] did the similar thing over the derivation Lie algebra
of the algebra of quantum torus.  The second author [Z] determined
the module structure of Shen's mixed product over Xu's nongraded Lie
algebras of Witt type in [X].

Let $M$ be any $o(n,\mbb{C})$-module and let $b\in\mbb{C}$.  It can
be verified that $\Im(\vt(o(n+2,\mbb{C})))\subset
\vt(o(n+2,\mbb{C}))+o(n,{\cal A})+{\cal A}\sum_{i=1}^nE_{i,i}$,
which enable us to define an $o(n+2,\mbb{C})$-module structure on
$\widehat M={\cal A}\otimes_\mbb{C}M$, where the hidden central
transformation $(\sum_{i=1}^nE_{i,i})|_{\widehat M}$ takes the
constant value $b$. We call such a module $\widehat M$ a {\it
generalized conformal module of} $o(n+2,\mbb{C})$. Geometrically,
the tensor modules yield natural
  extensions of the conformal transformations on
 a Riemannian manifold to its vector bundles. For any two integers $m$ and $n$, we use the notion
$$\overline{m,n}=\left\{\begin{array}{ll}\{m,m+1,\cdots,n\}&\mbox{if}\;m\leq n,\\
\emptyset&\mbox{otherwise}.\end{array}\right.\eqno(1.5)$$ Moreover,
we denote by $\mbb{N}$ the set of nonnegative integers.

Suppose that ${\cal E}$ is a vector space over $\mbb{R}$ with a
basis $\{\ves_i\mid i\in\ol{1,m}\}$ and the inner product
$(\sum_{i=1}^ma_i\ves_i,\sum_{r=1}^mc_r\ves_r)=\sum_{i=1}^na_ic_i$.
As usual, we take the set of simple positive roots
$\Pi_{2m}=\{\ves_i-\ves_{i+1},\ves_{m-1}+\ves_m\mid
i\in\ol{1,m-1}\}$ of $o(2m,\mbb{C})$ and the set of simple positive
roots $\Pi_{2m+1}=\{\ves_i-\ves_{i+1},\ves_m\mid i\in\ol{1,m-1}\}$
of $o(2m+1,\mbb{C})$. Set
$$\Lmd^+_n=\{\mu\in {\cal
E}\mid
[2(\mu,\nu)/(\nu,\nu)]\in\mbb{N}\;\for\;\nu\in\Pi_n\},\;\;n=2m,2m+1.\eqno(1.6)$$
It is well known that any finite-dimensional irreducible
$o(n,\mbb{C})$-module is a highest-weight irreducible module
$V(\mu)$ with highest weight $\mu\in\Lmd^+_n$ (e.g., cf. [Hu]). For
$\mu=\sum_{i=1}^m\mu_i\ves_i \in\Lmd^+_n$, we define
$\iota_\mu\in\ol{1,m}$ by
$$\mu_1=\mu_2=\cdots=\mu_{\iota_\mu}\neq
\mu_{\iota_\mu+1}.\eqno(1.7)$$ Using  Pieri's formulas, invariant
operators and the idea of Kostant's characteristic identities, we
prove: \psp

{\bf Main Theorem}. {\it Let $0\neq \mu=\sum_{i=1}^m\mu_i\ves_i
\in\Lmd^+_n$ with $n=2m,2m+1\geq 3$. The generalized conformal
module $\widehat{V(\mu)}$ of $o(2m,\mbb{C})$ is irreducible if $b\in
\mbb{C}\setminus\{m-1-\mbb{N}/2,\mu_1+2m-\iota_\mu-1-\mbb{N}\}$.
Moreover, the generalized conformal module $\widehat{V(\mu)}$ of
$o(2m+1,\mbb{C})$ is irreducible if $b\in
\mbb{C}\setminus\{m-\mbb{N}/2,\mu_1+2m-\iota_\mu-\mbb{N}\}$.

The generalized conformal module $\widehat{V(0)}$ of $o(n,\mbb{C})$
is irreducible if and only if $b\not\in -\mbb{N}$. When $b=0$,
$\widehat{V(0)}$ is isomorphic to the natural conformal
$o(n+2,\mbb{C})$-module ${\cal A}$, on which $A(f)=\vt(A)(f)$ for
$A\in o(n+2,\mbb{C})$ and $f\in{\cal A}$. The subspace $\mbb{C}$
forms a trivial $o(n+2,\mbb{C})$-submodule of the conformal module
${\cal A}$ and the quotient space ${\cal A}/\mbb{C}$ forms an
irreducible $o(n+2,\mbb{C})$-module.} \psp

In some singular cases of $\mu\neq 0$, the condition can be relaxed
slightly. We speculate that the module $\widehat{V(0)}$ may play the
same role in  higher-dimensional conformal field theory as that of
the vacuum module of the Virasoro algebra plays in two-dimensional
conformal field theory. The central charge in two-dimensional
conformal field theory may be replaced by the constant $b$ in the
above theorem for  higher-dimensional conformal field theory.

The paper is organized as follows. In Section 2, we prove the main
theorem for $o(2m,\mbb{C})$ ($m$ will be replaced by $n$
customarily). Section 3 is devoted to the proof of the main theorem
for $o(2m+1,\mbb{C})$ ($m$ will also be replaced by $n$
customarily).

\section{Generalized Conformal Representations of $D_{n+1}$ }

\quad\quad

Let $n>1$ be an integer.  The orthogonal Lie algebra
\begin{eqnarray*}\qquad o(2n,\mbb{C})&=&\sum_{1\leq
p<q\leq n}[\mbb{C}(E_{p,n+q}-E_{q,n+p})+
\mbb{C}(E_{n+p,q}-E_{n+q,p})]\\ &
&+\sum_{i,j=1}^n\mbb{C}(E_{i,j}-E_{n+j,n+i}).\hspace{7.3cm}(2.1)\end{eqnarray*}
We take the subspace
$${\cal H}=\sum_{i=1}^n\mbb{C}(E_{i,i}-E_{n+i,n+i})\eqno(2.2)$$
as a Cartan subalgebra and define $\{\ves_i\mid
i\in\ol{1,n}\}\subset{\cal H}^\ast$ by
$$\ves_i(E_{j,j}-E_{n+j,n+j})=\dlt_{i,j}.\eqno(2.3)$$
The inner product $(\cdot,\cdot)$ on the $\mbb{Q}$-subspace
$$L_\mbb{Q}=\sum_{i=1}^n\mbb{Q}\ves_i\eqno(2.4)$$
is given by
$$(\ves_i,\ves_j)=\dlt_{i,j}\qquad\for\;\;i,j\in\ol{1,n}.\eqno(2.5)$$
Then the root system of $o(2n,\mbb{C})$ is
$$\Phi_{D_n}=\{\pm \ves_i\pm\ves_j\mid1\leq i<j\leq
n\}.\eqno(2.6)$$ We take the set of positive roots
$$\Phi_{D_n}^+=\{\ves_i\pm\ves_j\mid1\leq i<j\leq
n\}.\eqno(2.7)$$ In particular,
$$\Pi_{D_n}=\{\ves_1-\ves_2,...,\ves_{n-1}-\ves_n,\ves_{n-1}+\ves_n\}\;\;\mbox{is the set of positive simple roots}.\eqno(2.8)$$

Recall the set of dominate integral weights
$$\Lmd^+=\{\mu\in L_\mbb{Q}\mid
(\ves_{n-1}+\ves_n,\mu),(\ves_i-\ves_{i+1},\mu)\in\mbb{N}\;\for\;i\in\ol{1,n-1}\}.\eqno(2.9)$$
According to (2.5),
$$\Lmd^+=\{\mu=\sum_{i=1}^n\mu_i\ves_i\mid
\mu_i\in\mbb{Z}/2;\mu_i-\mu_{i+1},\mu_{n-1}+\mu_n\in\mbb{N}\}.\eqno(2.10)$$
Note that if $\mu\in\Lmd^+$, then $\mu_{n-1}\geq |\mu_n|$. Given
$\mu\in\Lmd^+$, there exists a unique sequence ${\cal
S}(\mu)=\{n_0,n_1,...,n_s\}$ such that
$$0=n_0<n_1<n_2<\cdots <n_{s-1}<n_s=n\eqno(2.11)$$
and
$$\mu_i=\mu_j\Leftrightarrow n_r<i,j\leq n_{r+1}\;\;\mbox{for
some}\;\;r\in\ol{0,s-1}.\eqno(2.12)$$

For $\lmd\in\Lmd^+$, we denote by $V(\lmd)$ the finite-dimensional
irreducible $o(2n,\mbb{C})$-module with highest weight $\lmd$. The
$2n$-dimensional natural module of $o(2n,\mbb{C})$ is $V(\ves_1)$
with the weights $\{\pm\ves_i\mid i\in\ol{1,n}\}.$ The following
result is well known (e.g., cf. [FH]):\psp

{\bf Lemma 2.1 (Pieri's formula)}. {\it Given $\mu\in\Lmd^+$ with
${\cal S}(\mu)=\{n_0,n_1,...,n_s\}$,
$$V(\ves_1)\otimes_\mbb{C}V(\mu)\cong \bigoplus_{i=1}^s
(V(\mu-\ves_{n_i})\oplus V(\mu+\ves_{1+n_{i-1}}))\eqno(2.13)$$ if
$\mu_{n-1}+\mu_n>0$ and
$$V(\ves_1)\otimes_\mbb{C}V(\mu)\cong \bigoplus_{i=1}^{s-2+\dlt_{\mu_n,0}}
V(\mu-\ves_{n_i})\oplus \bigoplus_{i=1}^s
V(\mu+\ves_{1+n_{i-1}})\eqno(2.14)$$ when $\mu_{n-1}+\mu_n=0$.}\psp

Denote by $U({\cal G})$ the universal enveloping algebra of a Lie
algebra ${\cal G}$. The algebra $U({\cal G})$  can be imbedded into
the tensor algebra $U({\cal G})\otimes U({\cal G})$ by the
associative algebra homomorphism $d: U({\cal G}) \rightarrow U({\cal
G})\otimes_\mbb{C} U({\cal G})$ determined  by
$$d(u)=u\otimes 1 +1 \otimes u \qquad \mbox{ for} \ u\in
\cal{G}.\eqno(2.15)$$

Note that the Casimir element of $o(2n,\mbb{C})$ is
\begin{eqnarray*}\omega&=&\sum_{1\leq
i<j\leq
n}[(E_{i,n+j}-E_{j,n+i})(E_{n+j,i}-E_{n+i,j})+(E_{n+j,i}-E_{n+i,j})(E_{i,n+j}-E_{j,n+i})]
\\ & &+\sum_{i,j=1}^n(E_{i,j}-E_{n+j,n+i})(E_{j,i}-E_{n+i,n+j})\in
U(o(2n,\mbb{C})).\hspace{3.9cm}(2.16)\end{eqnarray*} Set
$$\td\omega=\frac{1}{2}(d(\omega)-\omega\otimes 1-1\otimes
\omega)\in U(o(2n,\mbb{C}))\otimes_\mbb{C}
U(o(2n,\mbb{C})).\eqno(2.17)$$ By (2.16),
\begin{eqnarray*}\td\omega\!&=&\!\sum_{1\leq
i<j\leq
n}[(E_{i,n+j}-E_{j,n+i})\otimes(E_{n+j,i}-E_{n+i,j})+(E_{n+j,i}-E_{n+i,j})\otimes(E_{i,n+j}-E_{j,n+i})]
\\ & &+\sum_{i,j=1}^n(E_{i,j}-E_{n+j,n+i})\otimes(E_{j,i}-E_{n+i,n+j}).\hspace{6.1cm}(2.18)\end{eqnarray*}
Denote
$$\rho=\frac{1}{2}\sum_{\nu\in\Phi_{D_n}^+}\nu.\eqno(2.19)$$
Then
$$(\rho,\nu)=1\qquad\for\;\;\nu\in\Pi_{D_n}\eqno(2.20)$$
(e.g., cf. [Hu]). By (2.8),
$$\rho=\sum_{i=1}^{n-1}(n-i)\ves_i.\eqno(2.21)$$

For any $\mu\in \Lmd^+$, we have
$$\omega|_{V(\mu)}=(\mu+2\rho,\mu)\mbox{Id}_{V(\mu)}.\eqno(2.22)$$
Denote
$$\ell^+_i(\mu)=\dim
V(\mu+\ves_i)\qquad\mbox{if}\;\;\mu+\ves_i\in\Lmd^+\eqno(2.23)$$ and
$$\ell^-_i(\mu)=\dim
V(\mu-\ves_i)\qquad\mbox{if}\;\;\mu-\ves_i\in\Lmd^+.\eqno(2.24)$$
Observe that
$$(\mu+\ves_i+2\rho,\mu+\ves_i)-(\mu+2\rho,\mu)-(\ves_1+2\rho,\ves_1)=2(\mu_i+1-i)\eqno(2.25)$$
and
$$(\mu-\ves_i+2\rho,\mu-\ves_i)-(\mu+2\rho,\mu)-(\ves_1+2\rho,\ves_1)=2(1+i-2n-\mu_i)\eqno(2.26)$$
for $\mu=\sum_{r=1}^n\mu_r\ves_r$ by (2.21). Moreover, the algebra
$U(o(2n,\mbb{C}))\otimes_\mbb{C}U(o(2n,\mbb{C}))$ acts on
$V(\ves_1)\otimes_\mbb{C}V(\mu)$ by
$$(\xi_1\otimes \xi_2)(v\otimes
u)=\xi_1(v)\otimes\xi_2(u)\qquad\for\;\;\xi_1,\xi_2\in
U(o(2n,\mbb{C})),\;v\in V(\ves_1),\;u\in V(\mu).\eqno(2.27)$$ By
Lemma 2.1, (2.17) and (2.23)-(2.26), we obtain:\psp

{\bf Lemma 2.2}. {\it Let $\mu=\sum_{i=1}^n\mu_i\ves_i\in\Lmd^+$
with ${\cal S}(\mu)=\{n_0,n_1,...,n_s\}$. If $\mu_{n-1}+\mu_n>0$,
the characteristic polynomial of
$\td\omega|_{V(\ves_1)\otimes_\mbb{C}V(\mu)}$ is
$$\prod_{i=1}^s(t-\mu_{1+n_{i-1}}+n_{i-1})^{\ell^+_{1+n_{i-1}}(\mu)}
(t+\mu_{n_i}+2n-n_i-1)^{\ell^-_{n_i}(\mu)}.\eqno(2.28)$$ When
$\mu_{n-1}+\mu_n=0$, the characteristic polynomial of
$\td\omega|_{V(\ves_1)\otimes_\mbb{C}V(\mu)}$ is
$$[\prod_{i=0}^{s-1}(t-\mu_{1+n_i}+n_i)^{\ell^+_{1+n_{i-1}}(\mu)}]
[\prod_{j=1}^{s-2+\dlt_{\mu_n,0}}(t+\mu_{n_j}+2n-n_j-1)^{\ell^-_{n_j}(\mu)}].\eqno(2.29)$$}
\pse

We remark that the above lemma is equivalent to special detailed
version of Kostant's characteristic identity. Set
$${\cal A}=\mbb{C}[x_1,x_2,...,x_{2n}].\eqno(2.30)$$
Then ${\cal A}$ forms an $o(2n,\mbb{C})$-module with the action
determined via
$$E_{i,j}|_{\cal
A}=x_i\ptl_{x_j}\qquad\for\;\;i,j\in\ol{1,2n}.\eqno(2.31)$$ The
corresponding Laplace operator and dual invariant are
$$\Dlt=\sum_{r=1}^n\partial_{x_r}\partial_{x_{n+r}},\qquad\eta=\sum_{i=1}^nx_ix_{n+i}.\eqno(2.32)$$Denote
$$D=\sum_{r=1}^{2n}x_r\partial_{x_r},\;\;J_i=x_iD-\eta\partial_{x_{n+i}},
\;\;J_{n+i}= x_{n+i}D-\eta\partial_{x_i}, \eqno(2.33)$$
$$A_{i,j}=x_i\partial_{x_j}-x_{n+j}\partial_{x_{n+i}} ,\;\;
B_{i,j}=x_i\partial_{x_{n+j}}-x_j\partial_{x_{n+i}} ,
C_{i,j}=x_{n+i}\partial_{x_j}-x_{n+j}\partial_{x_i}\eqno(2.34)
$$ for $i,j\in\ol{1,n}$. Then
$${\cal C}_{2n}=\mbb{C}D+\sum_{r=1}^{2n}(\mbb{C}\ptl_{x_r}+\mbb{C}J_r)+\sum_{i,j=1}^n\mbb{C}A_{i,j}+\sum_{1\leq
i<j\leq n}(\mbb{C}B_{i,j}+\mbb{C}C_{i,j})\eqno(2.35)$$ is the Lie
algebra of $2n$-dimensional  conformal group over $\mbb{C}$ with
$\eta$ as the metric.

Set
$${\cal L}_0=\sum_{i,j=1}^n\mathbb{C}A_{i,j}+\sum_{1\leq i<j\leq
n}(\mbb{C}B_{i,j}+\mathbb{C}C_{i,j}),\;\; {\cal
J}=\sum_{i=1}^{2n}\mathbb{C}J_{i}, \;\;{\cal
D}=\sum_{i=1}^{2n}\mathbb{C}\partial_{x_{i}}. \eqno(2.36)$$ Then
${\cal L}_0=o(2n,\mbb{C})|_{\cal A}$. We can easily verify  the
following Lie brackets:
$$[{\cal J},{\cal J}]=\{0\},\;\; [{\cal D}, {\cal D}]=\{0\},\
[\partial_{x_{k}},J_{n+i}]=C_{i,k},\;\;[\partial_{x_{n+k}},J_i]=B_{i,k},\eqno(2.37)$$
$$[\partial_{x_{k}},J_i]=\delta_{k,i}D+A_{i,k}, \qquad
[\partial_{x_{n+k}},J_{n+i}]=\delta_{k,i}D-A_{k,i}, \eqno(2.38)$$
$$[\partial_{x_{k}},A_{i,j}]=\delta_{k,i}\partial_{x_{j}},\;\;
 [\partial_{x_{k}},B_{i,j}]=\delta_{k,i}\partial_{x_{n+j}}-\delta_{k,j}
 \partial_{x_{n+i}},\;\;[\partial_{x_{k}},C_{i,j}]=0,\eqno(2.39)$$
 $$ [\partial_{x_{n+k}},A_{i,j}]=-\delta_{k,j}\partial_{x_{n+i}}, \;\; [\partial_{x_{n+k}},B_{i,j}]=0,
 \;\;[\partial_{x_{n+k}},C_{i,j}]=\delta_{k,i}\partial_{x_{j}}-\delta_{k,j}
 \partial_{x_{i}},\eqno(2.40)$$
$$[J_{k},A_{i,j}]=-\delta_{k,j}J_{i},\; \; [J_{k},B_{i,j}]=0, \;\; [J_{k},C_{i,j}]
=\delta_{k,i}J_{n+j}-\delta_{k,j}J_{n+i},\eqno(2.41) $$
$$[J_{n+k},A_{i,j}]=\delta_{k,i}J_{n+j}, \;\; [J_{n+k},B_{i,j}]=\delta_{k,i}J_j-\delta_{k,j}J_i, \;\;
  [J_{n+k},C_{i,j}]=0, \eqno(2.42)$$
$$[D,J_{k}]=J_{k}, \;\; [D, J_{n+k}]=J_{n+k}, \;\; [ \partial_{x_{k}},D]=\partial_{x_{k}}, \;\;
 [ \partial_{x_{n+k}},D]=\partial_{x_{n+k}}\eqno(2.43)$$
for $i,j,k\in\ol{1,n}$. Recall that the split \begin{eqnarray*}
\qquad o(2n+2,\mbb{C})&=&\sum_{1\leq r<s\leq
n+1}[\mbb{C}(E_{r,n+1+s}-E_{s,n+1+r})+\mbb{C}(E_{n+1+r,s}-E_{n+1+s,r})]\\
&&+
\sum_{i,j=1}^{n+1}\mbb{C}(E_{i,j}-E_{n+1+j,n+1+i}).\hspace{5.7cm}(2.44)\end{eqnarray*}
By (2.37)-(2.43), we have the  Lie algebra isomorphism
$\vt:o(2n+2,\mbb{F})\rightarrow {\cal C}_{2n}$ determined by
$$\vt(E_{i,j}-E_{n+1+j,n+1+i})=A_{i,j},\qquad\vt(E_{r,n+1+s}-E_{s,n+1+r})=B_{r,s},\eqno(2.45)$$
$$\vt(E_{n+1+r,s}-E_{n+1+s,r})=C_{r,s},\qquad \vt(E_{n+1,n+1}-E_{2n+2,2n+2})=-D,\eqno(2.46)$$
$$\vt(E_{n+1,i}-E_{n+1+i,2n+2})=\ptl_{x_i},\qquad\vt(E_{i,2n+2}-E_{n+1,n+1+i})=-\ptl_{x_{n+i}},\eqno(2.47)$$
$$\vt(E_{i,n+1}-E_{2n+2,i+n+1})=-J_i,
\qquad\vt(E_{2n+2,i}-E_{n+1+i,n+1})=J_{n+i}\eqno(2.48)$$ for
$i,j\in\ol{1,n}$ and $1\leq r<s\leq n$.

Recall the Witt algebra ${\cal W}_{2n}=\sum_{i=1}^{2n}{\cal
A}\ptl_{x_i}$, and Shen [Sg1-3] found a monomorphism $\Im$ from the
Lie algebra ${\cal W}_{2n}$ to the Lie algebra of semi-product
${\cal W}_{2n}+gl(2n,{\cal A})$ defined by
$$\Im(\sum_{i=1}^{2n}f_i\ptl_{x_i})=\sum_{i=1}^{2n}f_i\ptl_{x_i}+\sum_{i,j=1}^n\ptl_{x_i}(f_j)E_{i,j}.
\eqno(2.49)$$ Note that ${\cal C}_{2n}\subset {\cal W}_{2n}$ and
$$\Im(A_{i,j})=A_{i,j}+E_{i,j}-E_{n+j,n+i},
\Im(B_{r,s})=B_{r,s}+E_{r,n+s}-E_{s,n+r},\;\;\Im(\ptl_{x_i})=\ptl_{x_i},\eqno(2.50)$$
$$\Im(C_{r,s})=C_{r,s}+E_{n+r,s}-E_{n+s,r},\;\;\Im(\ptl_{x_{n+i}})=\ptl_{x_{n+i}},\;\;
\Im(D)=D+\sum_{p=1}^{2n}E_{p,p}\eqno(2.51)$$
$$\Im(J_{i})=J_{i}+\sum_{p=1}^nx_{n+p}(E_{i,n+p}-E_{p,n+i})+\sum_{q=1}^nx_q(E_{i,q}-E_{n+q,n+i})+x_i
\sum_{p=1}^{2n}E_{p,p},\eqno(2.52)$$
$$\Im(J_{n+i})=J_{n+i}+\sum_{p=1}^nx_{n+p}(E_{n+i,n+p}-E_{p,i})+\sum_{q=1}^nx_q(E_{n+i,q}-E_{n+q,i})+x_{n+i}
\sum_{p=1}^{2n}E_{p,p}\eqno(2.53)$$ for $i,j\in\ol{1,n}$ and $1\leq
r<s\leq n$. Moreover,
$$\widehat{\cal C}_{2n}={\cal C}_{2n}+o(2n,{\cal A})+{\cal
A}\sum_{p=1}^{2n}E_{p,p}\eqno(2.54)$$ forms a Lie subalgebra of
${\cal W}_{2n}+gl(2n,{\cal A})$ and $\Im({\cal C}_{2n})\subset
\widehat{\cal C}_{2n}$. In particular, the element
$\sum_{p=1}^{2n}E_{p,p}$ is a hidden central element.

Let  $M$ be an $o(2n,\mathbb{C})$-module and let $b\in\mbb{C}$ be a
fixed constant. Then
$$\widehat M={\cal A}\otimes_{\mbb{F}}M\eqno(2.55)$$
becomes a $\widehat{\cal C}_{2n}$-module with the action:
$$(d+f_1A+f_2\sum_{p=1}^{2n}E_{p,p})(g\otimes
v)=(d(g)+bf_2g)\otimes v+f_1g\otimes A(v)\eqno(2.56)$$ for
$f_1,f_2,g\in{\cal A},\;A\in o(2n,\mbb{C})$ and $v\in M$. Moreover,
we make $\widehat M$ a ${\cal C}_{2n}$-module with the action:
$$\xi(w)=\Im(\xi)(w)\qquad\for\;\;\xi\in {\cal C}_{2n},\;w\in
\widehat M.\eqno(2.57)$$ Furthermore, $\widehat M$ becomes an
$o(2n+2,\mbb{F})$-module with the action
$$A(w)=\Im(\vt(A))(w)\qquad\for\;\;A\in o(2n+2,\mbb{C}), \;w\in
\widehat M\eqno(2.58)$$ (cf. (2.45)-(2.48)). \psp

{\bf Lemma 2.3} \quad {\it If $M$ is an irreducible
$o(2n,\mbb{C})$-module, then the space $U({\cal J})(1\otimes M)$ is
an
 irreducible $o(2n+2,\mbb{C})$-submodule of $\widehat M$.}

{\it Proof}. Recall $o(2n,\mbb{C})|_{\cal A}={\cal L}_0$. Moreover,
(2.41) and (2.42) give
$$[{\cal L}_0,{\cal J}]={\cal J}.\eqno(2.59)$$
Note that $D$ is the degree operator (cf. (2.33)) and
$$[{\cal D},{\cal J}]\subset {\cal L}_0+\mbb{C}D\eqno(2.60)$$
by (2.37) and (2.38). According to (2.35) and (2.36),
$$\vt(o(2n+2,\mbb{C}))={\cal C}_{2n}={\cal L}_0+{\cal D}+{\cal
J}+\mbb{C}D.\eqno(2.61)$$ By (2.43) and (2.59)-(2.61), $U({\cal
J})(1\otimes M)$ forms an $o(2n+2,\mbb{C})$-submodule of $\widehat
M$.

Let $W$ be a nonzero $o(2n+2,\mbb{C})$-submodule of $U({\cal
J})(1\otimes M)$. Note
$$\ptl_{x_i}|_{\widehat M}=\ptl_{x_i}\otimes
1\qquad\for\;\;i\in\ol{1,2n}.\eqno(2.62)$$ By repeatedly applying
the above operators to $W$, we can prove
$$W\bigcap (1\otimes M)\neq \{0\}.\eqno(2.63)$$
However, $W\bigcap (1\otimes M)$ is a nonzero ${\cal L}_0$-submodule
of $1\otimes M$, which is an irreducible ${\cal L}_0$-module.
Therefore, $1\otimes M\subset W$. As a ${\cal C}_{2n}$-module,
$W\supset U({\cal J})(1\otimes M).\qquad\Box$\psp

 Write
$$x^{\alpha}=\prod_{i=1}^{2n}x_{i}^{\alpha_i} ,\;\;J^\al=\prod_{i=1}^{2n}J_{i}^{\alpha_i}
 \ \ \mbox{for} \
\al=(\al_1,\al_2,...,\al_{2n})\in\mbb{N}^{2n}.\eqno(2.64)$$
 For $k\in\mbb{N}$, we set
$${\cal A}_k=\mbox{Span}_{\mathbb{C}}\{x^\al\mid
\al \in\mbb{N}^{2n}, \ \sum_{i=1}^{2n}\al_i=k\},\;\; \widehat
M_{{\langle}k\rangle}={\cal A}_k\otimes_\mbb{C}M\eqno(2.65)$$ and
 $$(U({\cal J})(1\otimes
M))_{{\langle}k\rangle}=\mbox{Span}_{\mathbb{C}}\{ J^\al(1\otimes
M)\mid  \al \in\mbb{N}^{2n}, \ \sum_{i=1}^{2n}\al_i=k\}.
\eqno(2.67)$$
 Moreover,
$$(U({\cal J})(1\otimes
M))_{{\langle}0\rangle}=\widehat M_{{\langle}0\rangle}=1\otimes
M.\eqno(2.68)$$ Furthermore,
 $$\widehat M=\bigoplus\limits_{k=0}^\infty\widehat M_{\langle
 k\rangle},\qquad
 U({\cal J})(1\otimes M)=\bigoplus\limits_{k=0}^\infty(U({\cal J})(1\otimes
M))_{\langle k\rangle}.\eqno(2.69)$$

Next we define a linear transformation $\vf$ on  $\widehat M$
determined by
$$\vf(x^\al\otimes v)=J^\al(1\otimes
v)\qquad\for\;\;\al\in\mbb{N}^{2n},\;v\in M.\eqno(2.70)$$ Note that
${\cal A}_1=\sum_{i=1}^{2n}\mbb{C}x_i$ forms the $2n$-dimensional
natural ${\cal L}_0$-module (equivalently $o(2n,\mbb{C})$-module).
According to (2.41) and (2.42), ${\cal J}$ forms an ${\cal
L}_0$-module with respect to the adjoint representation, and the
linear map from ${\cal A}_1$ to ${\cal J}$ determined by $x_i\mapsto
J_i$ for $i\in\ol{1,2n}$ gives an ${\cal L}_0$-module isomorphism.
Thus $\vf$ can also be viewed as an ${\cal L}_0$-module homomorphism
from $\widehat M$ to $U({\cal J})(1\otimes M)$. Moreover,
$$\vf(\widehat M_{\la k\ra})=(U({\cal J})(1\otimes
M))_{{\langle}k\rangle}\qquad\for\;\;k\in\mbb{N}.\eqno(2.71)$$

{\bf Lemma 2.4}. {\it We have $\vf|_{\widehat M_{\la
1\ra}}=(b+\td\omega)|_{\widehat M_{\la 1\ra}}$ (cf.
(2.16)-(2.18)).}\psp

{\it Proof}.  Let $i\in \ol{1,n}$ and $v\in M$. Expressions (2.52),
(2.53) and (2.56)-(2.58) give
\begin{eqnarray*}\qquad\vf(x_i\otimes v)&=&J_i(1\otimes v)
=\sum_{p=1}^nx_{n+p}\otimes (E_{i,n+p}-E_{p,n+i})(v)\\
& &+\sum_{q=1}^nx_q\otimes (E_{i,q}-E_{n+q,n+i})(v)+bx_i\otimes v,
\hspace{3.9cm}(2.72)\end{eqnarray*}
\begin{eqnarray*}\qquad\vf(x_{n+i}\otimes v)&=&J_{n+i}(1\otimes v)=
\sum_{p=1}^nx_{n+p}\otimes(E_{n+i,n+p}-E_{p,i})(v)\\
& &+\sum_{q=1}^nx_q\otimes (E_{n+i,q}-E_{n+q,i})(v)+bx_{n+i}\otimes
v.\hspace{3.2cm}(2.73)\end{eqnarray*}

On the other hand, (2.18) and (2.31) yield
\begin{eqnarray*}& &\td\omega(x_i\otimes v)=
\sum_{1\leq p<q\leq
n}[(E_{p,n+q}-E_{q,n+p})\otimes(E_{n+q,p}-E_{n+p,q})(x_i\otimes v)
\\ & &+(E_{n+q,p}-E_{n+p,q})\otimes(E_{p,n+q}-E_{q,n+p})(x_i\otimes v)]
\\ & &+\sum_{r,s=1}^n(E_{r,s}-E_{n+s,n+r})\otimes(E_{s,r}-E_{n+r,n+s})(x_i\otimes v)
\\ &=&\sum_{p=1}^nx_{n+p}\otimes(E_{i,n+p}-E_{p,n+i})(v)
+\sum_{r=1}^nx_r\otimes(E_{i,r}-E_{n+r,n+i})(v),
\hspace{2.6cm}(2.74)\end{eqnarray*}
\begin{eqnarray*}& &\td\omega(x_{n+i}\otimes v)=
\sum_{1\leq p<q\leq
n}[(E_{p,n+q}-E_{q,n+p})\otimes(E_{n+q,p}-E_{n+p,q})(x_{n+i}\otimes
v)
\\ & &+(E_{n+q,p}-E_{n+p,q})\otimes(E_{p,n+q}-E_{q,n+p})(x_{n+i}\otimes v)]
\\ & &+\sum_{r,s=1}^n(E_{r,s}-E_{n+s,n+r})\otimes(E_{s,r}-E_{n+r,n+s})(x_{n+i}\otimes v)
\\ &=&\sum_{p=1}^nx_p\otimes(E_{n+i,p}-E_{n+p,i})(v)
+\sum_{s=1}^nx_{n+s}\otimes(E_{n+i,n+s}-E_{s,i})(v).
\hspace{2.7cm}(2.75)\end{eqnarray*} Comparing the above four
expressions, we get the conclusion in the lemma. $\qquad\Box$\psp

For $f\in{\cal A}$, we define the action
$$f(g\otimes v)=fg\otimes v\qquad\for\;\;g\in{\cal A},\;v\in
M.\eqno(2.76)$$ Then we have the $o(2n,\mbb{C})$-invariant operator
$$T=\sum_{i=1}^n(J_ix_{n+i}+J_{n+i}x_i)|_{\widehat M}.\eqno(2.77)$$

{\bf Lemma 2.5}. {\it We have $T|_{\widehat{M}_{\la
k\ra}}=(2b+2-2n+k)\eta|_{\widehat{M}_{\la k\ra}}$}.

{\it Proof}. Let $f$ be a homogeneous polynomial with degree $k$ and
let $v\in M$. By  (2.52), (2.53) and (2.57), we have
\begin{eqnarray*}& &T(f\otimes v)=\sum_{i=1}^n(J_ix_{n+i}+J_{n+i}x_i)(f\otimes v)
=\sum_{i=1}^n[J_i(x_{n+i}f\otimes v)+J_{n+i}(x_if\otimes v)]
\\ &=&\sum_{i=1}^n[J_i(x_{n+i}f)\otimes v+\sum_{p=1}^n(x_{n+p}x_{n+i}f)\otimes
(E_{i,n+p}-E_{p,n+i})(v)
\\ &
&+\sum_{q=1}^n(x_qx_{n+i}f)\otimes(E_{i,q}-E_{n+q,n+i})(v)+(x_ix_{n+i}f)\otimes
(\sum_{p=1}^{2n}E_{p,p})(v)
\\ & &+J_{n+i}(x_if)\otimes v+\sum_{p=1}^n(x_{n+p}x_i)f\otimes(E_{n+i,n+p}-E_{p,i})(v)
\\ & &+\sum_{q=1}^n(x_qx_if)\otimes(E_{n+i,q}-E_{n+q,i})(v)+(x_{n+i}x_if)\otimes
(\sum_{p=1}^{2n}E_{p,p})(v) ]
\\&=&
\sum_{i,p=1}^n[(x_{n+p}x_{n+i}f)\otimes
E_{i,n+p}(v)-(x_{n+p}x_{n+i}f)\otimes
E_{p,n+i}(v)]+\sum_{i,q=1}^n(x_qx_{n+i}f)\\
&&\otimes(E_{i,q}-E_{n+q,n+i})(v)
+\sum_{i,p=1}^n(x_{n+p}x_i)f\otimes(E_{n+i,n+p}-E_{p,i})(v)
+\sum_{i,q=1}^n[(x_qx_if)\\ & &\otimes E_{n+i,q}(v)-(x_qx_if)\otimes
E_{n+q,i}(v)]+[2b\eta f+\sum_{i=1}^n[J_i(x_{n+i}f)+J_{n+i}(x_if)]]\otimes v\\
&=&[\sum_{i=1}^n[J_i(x_{n+i}f)+J_{n+i}(x_if)]+2b\eta f]\otimes
v.\hspace{6.7cm}(2.78)
\end{eqnarray*}
According to (2.33), we find
\begin{eqnarray*}\qquad& &\sum_{i=1}^n[J_i(x_{n+i}f)+J_{n+i}(x_if)]
\\ &=&\sum_{i=1}^n[(x_iD-\eta\partial_{x_{n+i}})(x_{n+i}f)+
(x_{n+i}D-\eta\partial_{x_i})(x_if)]
\\ &=&2(k+1)\eta f-2n\eta f-\eta D(f)=(k+2-2n)\eta f.\hspace{4.5cm}(2.79)
\end{eqnarray*}
So the lemma holds.$\qquad\Box$\psp

For $ 0\neq\mu=\sum_{i=1}^n\mu_i\ves_i\in\Lmd^+$ with ${\cal
S}(\mu)=\{n_0,n_1,...,n_s\}$, we define
$$\Theta(\mu)=\left\{\begin{array}{ll}
\mu_1+n-1-\mbb{N}&\mbox{if}\;\mu_{n-1}=-\mu_n>0\;\mbox{and}\;s=2,
\\\mu_1+2n-n_1-1-\mbb{N}&\mbox{otherwise}.\end{array}\right.\eqno(2.80)$$

{\bf Theorem 2.6}.  {\it For $0\neq\mu\in\Lmd^+$, the generalized
conformal $o(2n+2,\mbb{C})$-module $\widehat{V(\mu)}$ defined by
(2.45)-(2.58) is irreducible if $b\in
\mbb{C}\setminus\{n-1-\mbb{N}/2,\Theta(\mu)\}.$}

{\it Proof}. By Lemma 2.3, it is enough to prove that the
homomorphism $\vf$ defined in (2.70) satisfies
$\vf(\widehat{V(\mu)})=\widehat{V(\mu)}$.  According to (2.71), we
only need to prove
$$\vf(\widehat{V(\mu)}_{\la k\ra})=\widehat{V(\mu)}_{\la k\ra}\eqno(2.81)$$
for any $k\in\mbb{N}$. We will prove it by induction on $k$.

When $k=0$, (2.81) holds by the definition (2.70). Consider $k=1$.
Write $\mu=\sum_{i=1}^n\mu_i\ves_i\in\Lmd^+$ with ${\cal
S}(\mu)=\{n_0,n_1,...,n_s\}$.  According to Lemma 2.2 and Lemma 2.4
with $M=V(\mu)$, the eigenvalues of $\vf|_{\widehat{V(\mu)}_{\la
1\ra}}$ are
$$b+\mu_{1+n_{i-1}}-n_{i-1},\;b-\mu_{n_i}-2n+n_i+1\;\;\for\;\;i\in\ol{1,s}\;\;\mbox{if}\;\;\mu_{n-1}+\mu_n>0\eqno(2.82)$$
and
$$b+\mu_{1+n_{i-1}}-n_{i-1},\;b-\mu_{n_r}-2n+n_r+1\;\;\for\;\;i\in\ol{1,s},\;r\in\ol{1,s-2+\dlt_{m_n,0}}\eqno(2.83)$$
when $\mu_{n-1}+\mu_n=0.$ Recall that $\mu_r\in\mbb{Z}/2$ for
$r\in\ol{1,n}$,
$$\mu_{\iota+1}-\mu_\iota\in\mbb{N}\;\;\for\;\;\iota\in\ol{1,n-1},\qquad
\;\;\mu_{n-1}\geq |\mu_n| \eqno(2.84)$$ and (2.12) holds. So
$$-\mu_{1+n_{i-1}}+n_{i-1},\;\mu_{n_i}+2n-n_i-1\in
\mu_1+2n-n_1-1-\mbb{N}\;\;\for\;\;i\in\ol{1,s}.\eqno(2.85)$$ If
$b\not\in \mu_1+2n-n_1-1-\mbb{N}$, then all the eigenvalues of
$\vf|_{\widehat{V(\mu)}_{\la 1\ra}}$ are nonzero.  In the case
$\mu_n=-\mu_{n-1}>0$ and $s=2$, $\mu_1=\mu_{n-1}=-\mu_n$ and the
eigenvalues $\vf|_{\widehat{V(\mu)}_{\la 1\ra}}$ are $b+\mu_1$ and
$b+\mu_n-n+1=b-\mu_1-n+1$, which are not equal to 0 because of
$b\not\in \Theta(\mu)=\mu_1+n-1-\mbb{N}$. Thus (2.81) holds for
$k=1$.

Suppose that (2.81) holds for $k\leq \ell$ with $\ell\geq 1$.
Consider $k=\ell+1$. Note that
$$\vf(\widehat{V(\mu)}_{\la \ell+1\ra})=\sum_{i=1}^{2n}\vf(x_i\widehat{V(\mu)}_{\la \ell\ra})
=\sum_{i=1}^{2n}J_i[\vf(\widehat{V(\mu)}_{\la
\ell\ra})]=\sum_{i=1}^{2n}J_i(\widehat{V(\mu)}_{\la
\ell\ra})\eqno(2.86)$$ by the inductional assumption.  To prove
(2.81) with $k=\ell+1$ is equivalent to prove
$$\sum_{i=1}^{2n}J_i(\widehat{V(\mu)}_{\la
\ell\ra})=\widehat{V(\mu)}_{\la \ell+1\ra}.\eqno(2.87)$$

For any $u\in \widehat{V(\mu)}_{\la \ell-1\ra}$, Lemma 2.5 says that
$$\sum_{i=1}^n[J_i(x_{n+i}u)+J_{n+i}(x_iu)]=(2b+1-2n+\ell)\eta
u.\eqno(2.88)$$ Since $b\not\in n-1-\mbb{N}/2$, $2b+1-2n+\ell\neq 0$
and  (2.88) gives
$$\eta u\in \sum_{i=1}^{2n}J_i(\widehat{V(\mu)}_{\la
\ell\ra})\qquad\for\;\;u\in  \widehat{V(\mu)}_{\la
\ell-1\ra}.\eqno(2.89)$$

Let $g\otimes v\in \widehat{V(\mu)}_{\la \ell\ra}$. According to
(2.52)-(2.57) and Lemma 2.4,
\begin{eqnarray*} J_i(g\otimes v)&=&x_iD(g)\otimes
v-\eta\ptl_{x_{n+i}}(g)\otimes
v+\sum_{p=1}^nx_{n+p}g\otimes(E_{i,n+p}-E_{p,n+i})(v)
\\& &+\sum_{q=1}^nx_qg\otimes(E_{i,q}-E_{n+q,n+i})(v)+x_ig\otimes
(\sum_{p=1}^{2n}E_{p,p})(v)\\ &=&-\eta\ptl_{x_{n+i}}(g)\otimes v
+g[(\ell+b+\td\omega)(x_i\otimes
v)]\hspace{5.3cm}(2.90)\end{eqnarray*} and
\begin{eqnarray*} J_{n+i}(g\otimes v)&=&x_{n+i}D(g)\otimes
v-\eta\ptl_{x_i}(g)\otimes
v+\sum_{p=1}^nx_{n+p}g\otimes(E_{n+i,n+p}-E_{p,i})(v)\\
& &+\sum_{q=1}^nx_qg\otimes(E_{n+i,q}-E_{n+q,i})(v)+x_{n+i}g\otimes
(\sum_{p=1}^{2n}E_{p,p})(v)\\ &=&-\eta\ptl_{x_i}(g)\otimes v
+g[(\ell+b+\td\omega)(x_{n+i}\otimes
v)]\hspace{4.9cm}(2.91)\end{eqnarray*} for $i\in\ol{1,n}$. Since
$$\eta\ptl_{x_i}(g)\otimes v,\;\eta\ptl_{x_{n+i}}(g)\otimes v
\in\sum_{r=1}^{2n}J_r(\widehat{V(\mu)}_{\la
\ell\ra})\qquad\for\;\;i\in\ol{1,n}\eqno(2.92)$$ by (2.89),
Expressions (2.90) and (2.91) show
$$g[(\ell+b+\td\omega)(x_i\otimes
v)]\in \sum_{r=1}^{2n}J_r(\widehat{V(\mu)}_{\la
\ell\ra})\qquad\for\;\;i\in\ol{1,2n},\;g\in{\cal
A}_\ell.\eqno(2.93)$$

According to Lemma 2.2 and 2.4, the eigenvalue of
$(\ell+b+\td\omega)|_{\widehat{V(\mu)}_{\la 1\ra}}$ are among
$$\{b+\ell+\mu_{1+n_{i-1}}-n_{i-1},\;b+\ell-\mu_{n_i}-2n+n_i+1\mid i\in\ol{1,s}\}.\eqno(2.94)$$
Again $$-\ell-\mu_{1+n_{i-1}}+n_{i-1},\;-\ell-\mu_{n_i}+2n-n_i-1\in
\mu_1+2n-n_1-1-\mbb{N}\;\;\for\;\;i\in\ol{1,s}.\eqno(2.95)$$ If
$b\not\in \mu_1+2n-n_1-1-\mbb{N}$, then all the eigenvalues of
$(\ell+b+\td\omega)|_{\widehat{V(\mu)}_{\la 1\ra}}$ are nonzero.  In
the case $\mu_n=-\mu_{n-1}>0$ and $s=2$, $\mu_1=\mu_{n-1}=-\mu_n$
and the eigenvalues $(\ell+b+\td\omega)|_{\widehat{V(\mu)}_{\la
1\ra}}$ are $b+\mu_1+\ell$ and $b+\mu_n-n+1+\ell=b-\mu_1-n+1+\ell$,
which are not equal to 0 because of $b\not\in
\Theta(\mu)=\mu_1+n-1-\mbb{N}$.
 Hence
$$(\ell+b+\td\omega)(\widehat{V(\mu)}_{\la 1\ra})=\widehat{V(\mu)}_{\la
1\ra}.\eqno(2.96)$$ By (2.93) and (2.96),
$$g(\widehat{V(\mu)}_{\la 1\ra})\subset \sum_{r=1}^{2n}J_r(\widehat{V(\mu)}_{\la
\ell\ra})\qquad\for\;\;g\in{\cal A}_\ell,\eqno(2.97)$$ equivalently,
(2.81) holds for $k=\ell+1$. By induction, (2.81) holds for any
$k\in\mbb{N}.\qquad\Box$\psp

We remark that the $o(2n+2,\mbb{C})$-module $\widehat{V(\mu)}$  is
$o(2n,\mbb{C})$-finite, that is, $\widehat{V(\mu)}$ is of $({\cal
G},{\cal K})$-type with ${\cal G}=o(2n+2,\mbb{C})$ and ${\cal
K}=o(2n,\mbb{C})$. Up to this stage, we do not known if the
condition in Theorem 2.6 is necessary for the generalized conformal
$o(2n+2,\mbb{C})$-module $\widehat{V(\mu)}$ to be irreducible if
$\mu\neq 0$. In the case $\mu=0$, the situation becomes clear. \psp

 {\bf Theorem 2.7}. {\it The generalized conformal
$o(2n+2,\mbb{C})$-module $\widehat{V(0)}$ is irreducible if and only
if $b\not\in-\mbb{N}$. When $b=0$, $\widehat{V(0)}$ is isomorphic to
the natural conformal $o(2n+2,\mbb{C})$-module ${\cal A}$, on which
$A(f)=\vt(A)(f)$ for $A\in o(2n+2,\mbb{C})$ and $f\in{\cal A}$ (cf.
(2.45)-(2.48)). The subspace $\mbb{C}$ forms a trivial
$o(2n+2,\mbb{C})$-submodule of the conformal module ${\cal A}$ and
the quotient space ${\cal A}/\mbb{C}$ forms an irreducible
$o(2n+2,\mbb{C})$-module.}

{\it Proof}. Pick $0\neq v_0\in V(0)$. Then $\widehat{V(0)}={\cal
A}\otimes v_0$.  Since $V(0)$ is the trivial $o(2n,\mbb{C})$-module,
(2.50)-(2.57) yield
$$\xi(f\otimes v_0)=\xi(f)\otimes v_0\qquad\for\;\;f\in {\cal A},\;\xi\in {\cal
L}_0+{\cal D}\eqno(2.98)$$ (cf. (2.36)) and
$$D(f\otimes v_0)=(b+D)(f)\otimes v_0,\eqno(2.99)$$
$$J_i(f\otimes v_0)=[x_i(D+b)-\eta\ptl_{x_{n+i}}](f)\otimes v_0,\;\;J_{n+i}
(f\otimes v_0)=[x_{n+i}(D+b)-\eta\ptl_{x_i}](f)\otimes
v_0\eqno(2.100)$$ for $f\in{\cal A}$ and $i\in\ol{1,2n}$ (cf.
(2.33)).

Recall the action of $o(2n+2,\mbb{C})$ on $\widehat{V(0)}$ by
(2.58). In particular, the map $f\otimes v_0\mapsto f$ for
$f\in{\cal A}$ gives an $o(2n,\mbb{C})$-module isomorphism from
$\widehat{V(0)}$ to ${\cal A}$. Remember the
$o(2n,\mbb{C})$-invariant differential operator
$\Dlt=\sum_{i=1}^n\ptl_{x_i}\ptl_{x_{n+i}}$ and its dual
$\eta=\sum_{i=1}^nx_ix_{n+i}$. Moreover, ${\cal A}_k$ denotes the
subspace of homogeneous polynomials in ${\cal A}$ with degree $k$.
Set
$${\cal H}_k=\{f\in{\cal
A}_k\mid\Dlt(f)=0\}\qquad\for\;\;k\in\mbb{N}.\eqno(2.101)$$ Then
${\cal H}_k\otimes v_0$ is an irreducible $o(2n,\mbb{C})$-submodule
with the highest-weight vector $x_1^k\otimes v_0$. Indeed,
$$\widehat{V(0)}=\bigoplus_{m,k=0}^\infty \eta^m{\cal H}_k\otimes
v_0\eqno(2.102)$$ is a direct sum of irreducible
$o(2n,\mbb{C})$-submodules. On the other hand, $U({\cal J})(1\otimes
v_0)$ forms an irreducible $o(2n+2,\mbb{C})$-submodule of
$\widehat{V(0)}$ (cf. Lemma 2.3). By (2.100), a necessary condition
for $\widehat{V(0)}=U({\cal J})(1\otimes v_0)$ is
$b\not\in-\mbb{N}$.

Next we assume $b\not\in -\mbb{N}-1$.  Let $W$ be a nonzero
$o(2n+2,\mbb{C})$-submodule of $\widehat{V(0)}$ such that
$$W\not\subset
\mbb{C}\otimes v_0\;\;\mbox{if}\;\; b= 0.\eqno(2.103)$$ By
repeatedly acting ${\cal D}$ on $W$ if necessary, we have $1\otimes
v_0\in W$. Note
$$J_i(1\otimes v_0)=bx_i\otimes
v_0\qquad\for\;\;i\in\ol{1,2n}.\eqno(2.104)$$ Thus
$$\widehat{V(0)}_{\la 1\ra}\subset W\eqno(2.105)$$
if $b\neq 0$. When $b=0$, (2.105) also holds because of (2.103),
$D(W)\subset W$ and the irreducibility of $\widehat{V(0)}_{\la
1\ra}$ as an $o(2n,\mbb{C})$-submodule. Suppose that
$$\widehat{V(0)}_{\la k\ra}\subset W\qquad\for\;\;k<\ell,\eqno(2.106)$$
where $2\leq \ell\in\mbb{N}$. According to (2.102),
$$\widehat{V(0)}_{\la \ell\ra}=\bigoplus_{m=0}^{\llbracket \ell/2
\rrbracket}\eta^m{\cal H}_{\ell-2m}\otimes v_0.\eqno(2.107)$$
Moreover,
$$[\Dlt,\eta]=n+D.\eqno(2.108)$$Set
$$\widehat{V(0)}_{\la \ell,r\ra}=\bigoplus_{m=0}^r\eta^m{\cal H}_{\ell-2m}\otimes v_0\eqno(2.109)$$
for $r\in\ol{0,\llbracket \ell/2 \rrbracket}$.  Then
$$\widehat{V(0)}_{\la \ell,r\ra}=\{w\in \widehat{V(0)}_{\la
\ell\ra}\mid \Dlt^{r+1}(w)=0\}\eqno(2.110)$$ and
$$\Dlt^r(\widehat{V(0)}_{\la \ell,r\ra})={\cal H}_{\ell-2r}\otimes
v_0\eqno(2.111)$$ by (2.108).

Since
$$J_1(x_1^{\ell-1}\otimes v_0)=(b+\ell-1)x_1^\ell\otimes v_0\in
W\eqno(2.112)$$ and $b\not\in-\mbb{N}-1$, we have
$$x_1^\ell\otimes v_0\in W.\eqno(2.113)$$
Hence
$$\widehat{V(0)}_{\la \ell,0\ra}={\cal H}_\ell\otimes v_0\subset
W\eqno(2.114)$$ because ${\cal H}_\ell\otimes v_0$ is an irreducible
$o(2n,\mbb{C})$-submodule generated by $x_1^\ell\otimes v_0$. Recall
$n\geq 2$ by our assumption. For $r\in\ol{1,\llbracket \ell/2
\rrbracket}$,
$$J_2(x_1^{\ell-r-1}x_{n+1}^r\otimes v_0)=(b+\ell-1)x_1^{\ell-r-1}x_2x_{n+1}^r\otimes v_0
\in W.\eqno(2.115)$$  So $x_2x_{n+1}^r\otimes v_0 \in W$. Moreover,
$$\Dlt^r(x_1^{\ell-r-1}x_2x_{n+1}^r\otimes
v_0)=r![\prod_{s=1}^r(\ell-r-s)]x_1^{\ell-2r-1}x_2\otimes v_0\in
{\cal H}_{\ell-2r}.\eqno(2.116)$$ Observe that (2.107) is a direct
sum of irreducible $o(2n,\mbb{C})$-submodules with distinct highest
weights. So
$$\widehat{V(0)}_{\la \ell\ra}\bigcap W=\bigoplus_{m=0}^{\llbracket \ell/2
\rrbracket}(\eta^m{\cal H}_{\ell-2m}\otimes v_0)\bigcap
W.\eqno(2.117)$$ By (2.109)-(2.111) and (2.116), $(\eta^r{\cal
H}_{\ell-2r}\otimes v_0)\bigcap W$ is a nonzero
$o(2n,\mbb{C})$-submodule. Since $\eta^r{\cal H}_{\ell-2r}\otimes
v_0$ is an irreducible $o(2n,\mbb{C})$-module, we have $\eta^r{\cal
H}_{\ell-2r}\otimes v_0=(\eta^r{\cal H}_{\ell-2r}\otimes v_0)\bigcap
W$. Therefore, $\widehat{V(0)}_{\la \ell\ra}\subset W$. By
induction, $\widehat{V(0)}_{\la k\ra}\subset W$ for any
$k\in\mbb{N}$, equivalently, $W=\widehat{V(0)}$. This proves the
theorem. $\qquad\Box$

\section{Generalized Conformal Representations of $B_{n+1}$ }

Let $n\geq 1$ be an integer.  The orthogonal Lie algebra
$$o(2n+1,\mbb{C})=o(2n,\mbb{C})+\sum_{i=1}^n[\mbb{C}(E_{0,i}-E_{n+i,0})
+\mbb{C}(E_{0,n+i}-E_{i,0})]\eqno(3.1)$$ (cf. (2.1)).  We take (2.2)
as a Cartan subalgebra and use the settings in (2.3)-(2.5). Then the
root system of $o(2n+1,\mbb{C})$ is
$$\Phi_{B_n}=\{\pm \ves_i\pm\ves_j,\pm\ves_r\mid1\leq i<j\leq
n;r\in\ol{1,n}\}.\eqno(3.2)$$ We take the set of positive roots
$$\Phi_{B_n}^+=\{\ves_i\pm\ves_j,\ves_r\mid1\leq i<j\leq
n,r\in\ol{1,n}\}.\eqno(3.3)$$ In particular,
$$\Pi_{B_n}=\{\ves_1-\ves_2,...,\ves_{n-1}-\ves_n,\ves_n\}\;\;\mbox{is the set of positive simple roots}.\eqno(3.4)$$

Recall the set of dominate integral weights
$$\Lmd^+=\{\mu\in L_\mbb{Q}\mid
2(\ves_n,\mu),(\ves_i-\ves_{i+1},\mu)\in\mbb{N}\;\for\;i\in\ol{1,n-1}.\eqno(3.5)$$
According to (2.5),
$$\Lmd^+=\{\mu=\sum_{i=1}^n\mu_i\ves_i\mid
\mu_i\in\mbb{N}/2;\mu_i-\mu_{i+1}\in\mbb{N}\;\for\;i\in\ol{1,n-1}\}.\eqno(3.6)$$
 Given
$\mu=\sum_{i=1}^n\mu_i\ves_i\in\Lmd^+$, there exists a unique
sequence ${\cal S}(\mu)=\{n_0,n_1,n_2,...,n_s\}$ such that (2.11)
and (2.12) holds. Denote
$$\rho=\frac{1}{2}\sum_{\nu\in\Phi_{B_n}^+}\nu.\eqno(3.7)$$
Then
$$\frac{2(\rho,\nu)}{(\nu,\nu)}=1\qquad\for\;\;\nu\in\Pi_{B_n}\eqno(3.8)$$
(e.g., cf. [Hu]). By (3.4),
$$\rho=\sum_{i=1}^{n-1}(n-i+1/2)\ves_i.\eqno(3.9)$$

 For $\lmd\in\Lmd^+$, we denote by $V(\lmd)$ the
finite-dimensional irreducible $o(2n+1,\mbb{C})$-module with highest
weight $\lmd$. The $(2n+1)$-dimensional natural module of
$o(2n+1,\mbb{C})$ is $V(\ves_1)$  with weights $\{0,\pm\ves_i\mid
i\in\ol{1,n}\}$. The following result is well known (e.g., cf.
[FH]):\psp

{\bf Lemma 3.1 (Pieri's formula)}. {\it Given $\mu\in\Lmd^+$ with
${\cal S}(\mu)=\{n_0,n_1,...,n_s\}$,
$$V(\ves_1)\otimes_\mbb{C}V(\mu)\cong
V(\mu)\oplus\bigoplus_{i=1}^{s-\dlt_{\mu_n,0}-\dlt_{\mu_n,1/2}}
V(\mu-\ves_{n_i})\oplus\bigoplus_{r=1}^s
V(\mu+\ves_{1+n_{r-1}}).\eqno(3.10)$$}\pse

Note that the Casimir element of $o(2n+1,\mbb{C})$ is
\begin{eqnarray*}\omega&=&\sum_{1\leq
i<j\leq
n}[(E_{i,n+j}-E_{j,n+i})(E_{n+j,i}-E_{n+i,j})+(E_{n+j,i}-E_{n+i,j})(E_{i,n+j}-E_{j,n+i})]
\\
&&+\sum_{i=1}^n[(E_{0,i}-E_{n+i,0})(E_{i,0}-E_{0,n+i})+(E_{i,0}-E_{0,n+i})(E_{0,i}-E_{n+i,0})]
\\ & &+\sum_{i,j=1}^n(E_{i,j}-E_{n+j,n+i})(E_{j,i}-E_{n+i,n+j})\in
U(o(2n+1,\mbb{C})).\hspace{3.2cm}(3.11)\end{eqnarray*} Set
$$\td\omega=\frac{1}{2}(d(\omega)-\omega\otimes 1-1\otimes
\omega)\in U(o(2n+1,\mbb{C}))\otimes_\mbb{C}
U(o(2n+1,\mbb{C})).\eqno(3.12)$$ By (3.11),
\begin{eqnarray*}\td\omega\!&=&\!\sum_{1\leq
i<j\leq
n}[(E_{i,n+j}-E_{j,n+i})\otimes(E_{n+j,i}-E_{n+i,j})+(E_{n+j,i}-E_{n+i,j})\otimes(E_{i,n+j}-E_{j,n+i})]
\\
&&+\sum_{i=1}^n[(E_{0,i}-E_{n+i,0})\otimes(E_{i,0}-E_{0,n+i})+(E_{i,0}-E_{0,n+i})\otimes(E_{0,i}-E_{n+i,0})]
\\ & &+\sum_{i,j=1}^n(E_{i,j}-E_{n+j,n+i})\otimes(E_{j,i}-E_{n+i,n+j}).\hspace{6.1cm}(3.13)\end{eqnarray*}
Moreover, (2.22) also holds. Take the settings in (2.23) and (2.24).
Denote
$$\ell(\mu)=\dim V(\mu).\eqno(3.14)$$

 Observe that
$$(\mu+2\rho,\mu)-(\mu+2\rho,\mu)-(\ves_1+2\rho,\ves_1)=-2n,\eqno(3.15)$$
$$(\mu+\ves_i+2\rho,\mu+\ves_i)-(\mu+2\rho,\mu)-(\ves_1+2\rho,\ves_1)=2(\mu_i+1-i)\eqno(3.16)$$
and
$$(\mu-\ves_i+2\rho,\mu-\ves_i)-(\mu+2\rho,\mu)-(\ves_1+2\rho,\ves_1)=2(i-2n-\mu_i)\eqno(3.17)$$
for $\mu=\sum_{r=1}^n\mu_r\ves_r$ by (3.9). Moreover, the algebra
$U(o(2n+1,\mbb{C}))\otimes_\mbb{C}U(o(2n+1,\mbb{C}))$ acts on
$V(\ves_1)\otimes_\mbb{C}V(\mu)$ by
$$(\xi_1\otimes \xi_2)(v\otimes
u)=\xi_1(v)\otimes\xi_2(u)\;\for\;\xi_1,\xi_2\in
U(o(2n+1,\mbb{C})),\;v\in V(\ves_1),\;u\in V(\mu).\eqno(3.18)$$ By
Lemma 3.1, (3.12), (2.22) and (3.15)-(3.17), we get:\psp

{\bf Lemma 3.2}. {\it Let $\mu=\sum_{i=1}^n\mu_i\ves_i\in\Lmd^+$
with ${\cal S}(\mu)=\{n_0,n_1,...,n_s\}$. The characteristic
polynomial of $\td\omega|_{V(\ves_1)\otimes_\mbb{C}V(\mu)}$ is
$$(t+n)^{\ell(\mu)}[\prod_{i=0}^{s-\dlt_{\mu_n,0}-\dlt_{\mu_n,1/2}}(t-\mu_{1+n_i}+n_i)^{\ell^+_{1+n_{i-1}}(\mu)}]
[\prod_{j=1}^{s-1}(t+\mu_{n_j}+2n-n_j-1)^{\ell^-_{n_j}(\mu)}].\eqno(3.19)$$}
\pse

We remark that the above lemma is also equivalent to special
detailed version of Kostant's characteristic identity. Set
$${\cal A}=\mbb{C}[x_0,x_1,x_2,...,x_{2n}].\eqno(3.20)$$
Then ${\cal A}$ forms an $o(2n+1,\mbb{C})$-module with the action
determined via
$$E_{i,j}|_{\cal
A}=x_i\ptl_{x_j}\qquad\for\;\;i,j\in\ol{0,2n}.\eqno(3.21)$$ The
corresponding Laplace operator and dual invariant are
$$\Dlt=\ptl_{x_0}^2+2\sum_{r=1}^n\partial_{x_r}\partial_{x_{n+r}},\qquad\eta=\frac{1}{2}x_0^2+\sum_{i=1}^nx_ix_{n+i}.
\eqno(3.22)$$Denote
$$D=\sum_{r=0}^{2n}x_r\partial_{x_r},\;\;J_i=x_iD-\eta\partial_{x_{n+i}},
\;\;J_{n+i}= x_{n+i}D-\eta\partial_{x_i}, \eqno(3.23)$$
$$J_0=x_0D-\eta\ptl_{x_0},\;\;K_i=x_0\ptl_{x_i}-x_{n+i}\ptl_{x_0},\;\;K_{n+i}=x_0\ptl_{x_{n+i}}-x_i\ptl_{x_0},
\eqno(3.24)$$
$$A_{i,j}=x_i\partial_{x_j}-x_{n+j}\partial_{x_{n+i}} ,\;\;
B_{i,j}=x_i\partial_{x_{n+j}}-x_j\partial_{x_{n+i}} ,
C_{i,j}=x_{n+i}\partial_{x_j}-x_{n+j}\partial_{x_i}\eqno(3.25)
$$ for $i,j\in\ol{1,n}$. Then
$${\cal C}_{2n+1}=\mbb{C}D+\sum_{s=1}^{2n}\mbb{C}K_s+\sum_{r=0}^{2n}(\mbb{C}\ptl_{x_r}+\mbb{C}J_r)+\sum_{i,j=1}^n\mbb{C}A_{i,j}+\sum_{1\leq
i<j\leq n}(\mbb{C}B_{i,j}+\mbb{C}C_{i,j})\eqno(3.26)$$ is the Lie
algebra of $(2n+1)$-dimensional  conformal group over $\mbb{C}$ with
$\eta$ in (3.22) as the metric. Moreover, the subspace
$${\cal C}'_{2n}=\mbb{C}D+\sum_{r=1}^{2n}(\mbb{C}\ptl_{x_r}+\mbb{C}J_r)+\sum_{i,j=1}^n\mbb{C}A_{i,j}+\sum_{1\leq
i<j\leq n}(\mbb{C}B_{i,j}+\mbb{C}C_{i,j})\eqno(3.27)$$ forms a Lie
subalgebra of ${\cal C}_{2n+1}$ which is isomorphic to ${\cal
C}_{2n}$ in (2.35).

 Set
$${\cal L}_0=\sum_{s=1}^{2n}\mbb{C}K_s+\sum_{i,j=1}^n\mathbb{C}A_{i,j}+\sum_{1\leq i<j\leq
n}(\mbb{C}B_{i,j}+\mathbb{C}C_{i,j}),\;\; {\cal
J}=\sum_{i=0}^{2n}\mathbb{C}J_{i}, \;\;{\cal
D}=\sum_{i=0}^{2n}\mathbb{C}\partial_{x_{i}}. \eqno(3.28)$$ Then
${\cal L}_0=o(2n+1,\mbb{C})|_{\cal A}$ and (2.37)-(2.43) hold.
Moreover,
$$[\ptl_{x_0},J_0]=D,\;\; [\ptl_{x_0}, J_{n+i}]=-K_i, \;\;
[\ptl_{x_0},J_i]=-K_{n+i},\;\;[\ptl_{x_{0}},K_{i}]=\ptl_{x_i},\eqno(3.29)$$
$$
[\ptl_{x_0},K_{n+i}]=\ptl_{x_{n+i}},\;\;[\ptl_{x_i},J_0]=K_i,\;\;[\ptl_{x_{n+i}},J_0]=K_{n+i},\eqno(3.30)$$
$$[\ptl_{x_i},K_j]=[\ptl_{x_{n+i}},K_{n+j}]=0,\;\;
[\ptl_{x_i},K_{n+j}]=[\ptl_{x_{n+i}},K_j]=-\dlt_{i,j}\ptl_{x_0},
\eqno(3.31)$$
$$[J_0,K_i]=J_{n+i},\;\;[J_0,K_{n+i}]=J_i,\;\;[J_i,K_j]=[J_{n+i},K_{n+j}]=-\dlt_{i,j}J_0,\eqno(3.32)$$
$$[J_i,K_{n+j}]=[J_{n+i},K_j]=0,\;\;[\ptl_{x_0},o(2n,\mbb{C})|_{\cal
A}]=[J_0,,o(2n,\mbb{C})|_{\cal A}]=\{0\}\eqno(3.33)$$ for
$i,j\in\ol{1,n}$ (cf. (3.21)).

Recall that the split
$$o(2n+3,\mbb{C})= o(2n+2,\mbb{C})+\sum_{i=1}^{n+1}[\mbb{C}(E_{0,i}-E_{n+i+1,0})+
(E_{0,n+1+i}-E_{i,0})]\eqno(3.34)$$ (cf. (2.44)). By (3.29)-(3.33),
we have the  Lie algebra isomorphism $\vt:o(2n+3,\mbb{F})\rightarrow
{\cal C}_{2n+1}$ determined by (2.45)-(2.48) and
$$\vt(E_{0,i}-E_{n+1+i,0})=K_i,\qquad
\vt(E_{0,n+1+i}-E_{i,0})=K_{n+i}\qquad\for\;i\in\ol{1,n},\eqno(3.35)$$
$$\vt(E_{0,2n+2}-E_{n+1,0})=-\ptl_{x_0},\qquad \vt(E_{0,n+1}-E_{2n+2,0})=J_0.\eqno(3.36)$$
Recall the Witt algebra ${\cal W}_{2n+1}=\sum_{i=0}^{2n}{\cal
A}\ptl_{x_i}$, and Shen [Sg1-3] found a monomorphism $\Im$ from the
Lie algebra ${\cal W}_{2n}$ to the Lie algebra of semi-product
${\cal W}_{2n+1}+gl(2n+1,{\cal A})$ defined by
$$\Im(\sum_{i=0}^{2n}f_i\ptl_{x_i})=\sum_{i=0}^{2n}f_i\ptl_{x_i}+\sum_{i,j=0}^n\ptl_{x_i}(f_j)E_{i,j}.
\eqno(3.37)$$ Note that ${\cal C}_{2n+1}\subset {\cal W}_{2n+1}$,
and (2.50) and (2.51) except the last equation hold. Moreover,
$$\Im(K_i)=K_i+E_{0,i}-E_{n+i,0},\qquad
\Im(D)=D+\sum_{p=0}^{2n}E_{p,p},\eqno(3.38)$$
$$\Im(K_{n+i})=K_{n+i}+E_{0,n+i}-E_{i,0},\qquad\Im(\ptl_{x_0})=\ptl_{x_0},\eqno(3.39)$$
$$\Im(J_0)=J_0+\sum_{s=1}^n[x_s(E_{0,s}-E_{n+s,0})+x_{n+s}(E_{0,n+s}-E_{s,0})]+x_0\sum_{p=0}^{2n}E_{p,p},\eqno(3.40)$$
\begin{eqnarray*}\Im(J_i)&=&J_{i}+\sum_{p=1}^nx_{n+p}(E_{i,n+p}-E_{p,n+i})+\sum_{q=1}^nx_q(E_{i,q}-E_{n+q,n+i})
\\ & &+x_0(E_{i,0}-E_{0,n+i})+x_i
\sum_{p=0}^{2n}E_{p,p},\hspace{7.1cm}(3.41)\end{eqnarray*}
\begin{eqnarray*}\Im(J_{n+i})&=&J_{n+i}+\sum_{p=1}^nx_{n+p}(E_{n+i,n+p}-E_{p,i})+\sum_{q=1}^nx_q(E_{n+i,q}-E_{n+q,i})
\\&&+x_0(E_{n+i,0}-E_{0,i})+x_{n+i} \sum_{p=0}^{2n}E_{p,p}.\hspace{6.3cm}(3.42)\end{eqnarray*}
Furthermore,
$$\widehat{\cal C}_{2n+1}={\cal C}_{2n+1}+o(2n+1,{\cal A})+{\cal
A}\sum_{p=0}^{2n}E_{p,p}\eqno(3.43)$$ forms a Lie subalgebra of
${\cal W}_{2n+1}+gl(2n+1,{\cal A})$ and $\Im({\cal C}_{2n+1})\subset
\widehat{\cal C}_{2n+1}$. In particular, the element
$\sum_{p=0}^{2n}E_{p,p}$ is a hidden central element.

Let  $M$ be an $o(2n+1,\mathbb{C})$-module and let $b\in\mbb{C}$ be
a fixed constant. Then
$$\widehat M={\cal A}\otimes_{\mbb{F}}M\eqno(3.44)$$
becomes a $\widehat{\cal C}_{2n+1}$-module with the action:
$$(d+f_1A+f_2\sum_{p=0}^{2n}E_{p,p})(g\otimes
v)=(d(g)+bf_2g)\otimes v+f_1g\otimes A(v)\eqno(3.45)$$ for
$f_1,f_2,g\in{\cal A},\;A\in o(2n+1,\mbb{C})$ and $v\in M$.
Moreover, we make $\widehat M$ a ${\cal C}_{2n}$-module with the
action:
$$\xi(w)=\Im(\xi)(w)\qquad\for\;\;\xi\in {\cal C}_{2n+1},\;w\in
\widehat M.\eqno(3.46)$$ Furthermore, $\widehat M$ becomes an
$o(2n+3,\mbb{F})$-module with the action
$$A(w)=\Im(\vt(A))(w)\qquad\for\;\;A\in o(2n+2,\mbb{C}), \;w\in
\widehat M\eqno(3.47)$$ (cf. (2.45)-(2.48), (3.35) and (3.36)). By a
proof similar to that of Lemma 2.3, we have:
 \psp

{\bf Lemma 3.3} \quad {\it If $M$ is an irreducible
$o(2n+1,\mbb{C})$-module, then the space $U({\cal J})(1\otimes M)$
is an
 irreducible $o(2n+3,\mbb{C})$-submodule of $\widehat M$.}
\psp

 Write
$$x^{\alpha}=\prod_{i=0}^{2n}x_{i}^{\alpha_i} ,\;\;J^\al=\prod_{i=0}^{2n}J_{i}^{\alpha_i}
 \ \ \mbox{for} \
\al=(\al_0,\al_1,...,\al_{2n})\in\mbb{N}^{2n+1}.\eqno(3.48)$$
 For $k\in\mbb{N}$, we set
$${\cal A}_k=\mbox{Span}_{\mathbb{C}}\{x^\al\mid
\al \in\mbb{N}^{2n+1}, \ \sum_{i=0}^{2n}\al_i=k\},\;\; \widehat
M_{{\langle}k\rangle}={\cal A}_k\otimes_\mbb{C}M\eqno(3.49)$$ and
 $$(U({\cal J})(1\otimes
M))_{{\langle}k\rangle}=\mbox{Span}_{\mathbb{C}}\{ J^\al(1\otimes
M)\mid  \al \in\mbb{N}^{2n+1}, \ \sum_{i=0}^{2n}\al_i=k\}.
\eqno(3.50)$$
 Moreover,
$$(U({\cal J})(1\otimes
M))_{{\langle}0\rangle}=\widehat M_{{\langle}0\rangle}=1\otimes
M.\eqno(3.51)$$ Furthermore,
 $$\widehat M=\bigoplus\limits_{k=0}^\infty\widehat M_{\langle
 k\rangle},\qquad
 U({\cal J})(1\otimes M)=\bigoplus\limits_{k=0}^\infty(U({\cal J})(1\otimes
M))_{\langle k\rangle}.\eqno(3.52)$$

Next we define a linear transformation $\vf$ on  $\widehat M$
determined by
$$\vf(x^\al\otimes v)=J^\al(1\otimes
v)\qquad\for\;\;\al\in\mbb{N}^{2n+1},\;v\in M.\eqno(3.53)$$ Note
${\cal A}_1=\sum_{i=0}^{2n}\mbb{C}x_i$  forms the
$(2n+1)$-dimensional natural ${\cal L}_0$-module (equivalently
$o(2n+1,\mbb{C})$-module). According to (2.41), (2.42), (3.32) and
(3.33), ${\cal J}$ forms an ${\cal L}_0$-module with respect to the
adjoint representation, and the linear map from ${\cal A}_1$ to
${\cal J}$ determined by $x_i\mapsto J_i$ for $i\in\ol{0,2n}$ gives
an ${\cal L}_0$-module isomorphism. Thus $\vf$ can also be viewed as
an ${\cal L}_0$-module homomorphism from $\widehat M$ to $U({\cal
J})(1\otimes M)$. Moreover,
$$\vf(\widehat M_{{\langle}k\rangle})=(U({\cal J})(1\otimes
M))_{{\langle}k\rangle}\qquad\for\;\;k\in\mbb{N}.\eqno(3.54)$$ \pse

{\bf Lemma 3.4}. {\it We have $\vf|_{\widehat M_{\la
1\ra}}=(b+\td\omega)|_{\widehat M_{\la 1\ra}}$ (cf.
(3.11)-(3.13)).}\psp

{\it Proof}.  Recall $\widehat M_{\la 1\ra}={\cal
A}_1\otimes_\mbb{C}M$. Let $i\in \ol{1,n}$ and $v\in M$. Expressions
(3.40), (3.45) and (3.46) give
$$\vf(x_0\otimes v)=\sum_{s=1}^n[x_s\otimes
(E_{0,s}-E_{n+s,0})(v)+x_{n+s}\otimes(E_{0,n+s}-E_{s,0})(v)]+bx_0\otimes
v.\eqno(3.55)$$ Moreover, (3.41), (3.45) and (3.46) imply
\begin{eqnarray*}\qquad\vf(x_i\otimes v)&=&\sum_{p=1}^nx_{n+p}\otimes (E_{i,n+p}-E_{p,n+i})(v)
+x_0\otimes(E_{i,0}-E_{0,n+i})(v)\\
& &+\sum_{q=1}^nx_q\otimes (E_{i,q}-E_{n+q,n+i})(v)+bx_i\otimes v
\hspace{4.1cm}(3.56)\end{eqnarray*}for $i\in\ol{1,n}$. Furthermore,
(3.42), (3.45) and (3.46) yield
\begin{eqnarray*}\qquad\vf(x_{n+i}\otimes v)&=&
\sum_{p=1}^nx_{n+p}\otimes(E_{n+i,n+p}-E_{p,i})(v)+x_0\otimes(E_{n+i,0}-E_{0,i})(v)\\
& &+\sum_{q=1}^nx_q\otimes (E_{n+i,q}-E_{n+q,i})(v)+bx_{n+i}\otimes
v\hspace{3.3cm}(3.57)\end{eqnarray*} for $i\in\ol{1,n}$.

On the other hand, (3.13) and (3.21) yield
$$\td\omega(x_0\otimes
v)=\sum_{i=1}^n[-x_{n+i}\otimes(E_{i,0}-E_{0,n+i})(v)
+x_i\otimes(E_{0,i}-E_{n+i,0})(v)],\eqno(3.58)$$
\begin{eqnarray*}\td\omega(x_i\otimes v)&=&
\sum_{p=1}^nx_{n+p}\otimes(E_{i,n+p}-E_{p,n+i})(v)+x_0\otimes(E_{i,0}-E_{0,n+i})(v)
\\ & &+\sum_{r=1}^nx_r\otimes(E_{i,r}-E_{n+r,n+i})(v),
\hspace{6.6cm}(3.59)\end{eqnarray*}
\begin{eqnarray*}\td\omega(x_{n+i}\otimes v)&=&
\sum_{p=1}^nx_p\otimes(E_{n+i,p}-E_{n+p,i})(v)-x_0\otimes(E_{0,i}-E_{n+i,0})(v)\\
& & +\sum_{s=1}^nx_{n+s}\otimes(E_{n+i,n+s}-E_{s,i})(v).
\hspace{5.85cm}(3.60)\end{eqnarray*}
 Comparing the above six
expressions, we get the conclusion in the lemma. $\qquad\Box$\psp

We use the definition (2.76) and  we have  the
$o(2n+1,\mbb{C})$-invariant operator
$$T=[J_0x_0+\sum_{i=1}^n(J_ix_{n+i}+J_{n+i}x_i)]|_{\widehat M}.\eqno(3.61)$$
\pse

{\bf Lemma 3.5}. {\it We have $T|_{\widehat{M}_{\la
k\ra}}=(2b-2n+k+1)\eta$}.

{\it Proof}. Let $f\in {\cal A}_k$ and  $v\in M$. According to
(3.24) and (3.40),
\begin{eqnarray*}J_0x_0(f\otimes v)&=&x_0\sum_{s=1}^n[x_sf\otimes(E_{0,s}-E_{n+s,0})(v)+x_{n+s}f\otimes
(E_{0,n+s}-E_{s,0})(v)]\\ & &+ [(k+1+b)x_0^2-\eta]f\otimes v-\eta
x_0\ptl_{x_0}(f)\otimes v.\hspace{4cm}(3.62)\end{eqnarray*}
Moreover, (2.79), (3.23), (3.41) and (3.42) give
\begin{eqnarray*}&&\sum_{i=1}^n(J_ix_{n+i}+J_{n+i}x_i)(f\otimes v)\\&=&
x_0\sum_{i=1}^n[x_if\otimes
(E_{n+i,0}-E_{0,i})(v)+x_{n+i}f\otimes(E_{i,0}-E_{0,n+i})(v)]\\ &&
+2[(b+k+1)(\sum_{i=1}^nx_ix_{n+i})-n\eta](f)\otimes
v-\eta(\sum_{i=1}^{2n}x_i\ptl_{x_i}(f)\otimes
v.\hspace{2.8cm}(3.63)\end{eqnarray*} Thus
$$T(f\otimes v)=[2(b+k+1)-2n-1]\eta(f)\otimes v-\eta D(f)\otimes
v=(2b+k+1-2n)\eta f\otimes v.\eqno(3.64)$$ So the Lemma
holds.$\qquad\Box$\psp

For $ 0\neq\mu=\sum_{i=1}^n\mu_i\ves_i\in\Lmd^+$ with ${\cal
S}(\mu)=\{n_0,n_1,...,n_s\}$, we define
$$\Theta(\mu)=\left\{\begin{array}{ll}\emptyset
&\mbox{if}\;\mu=(\sum_{i=1}^n\ves_i)/2,
\\\mu_1+2n-n_1-\mbb{N}&\mbox{otherwise}.\end{array}\right.\eqno(3.65)$$
\pse

 {\bf Theorem 3.6}.  {\it For $0\neq\mu\in\Lmd^+$, the
generalized conformal $o(2n+2,\mbb{C})$-module $\widehat{V(\mu)}$
defined by (2.50), (2.51) except the last equation,  and
(3.37)-(3.47) is irreducible if $b\in
\mbb{C}\setminus\{n-\mbb{N}/2,\Theta(\mu)\}.$}

{\it Proof}. By Lemma 3.3, it is enough to prove that the
homomorphism $\vf$ defined in (3.53) satisfies
$\vf(\widehat{V(\mu)})=\widehat{V(\mu)}$.  According to (3.54), we
only need to prove
$$\vf(\widehat{V(\mu)}_{\la k\ra})=\widehat{V(\mu)}_{\la k\ra}\eqno(3.66)$$
for any $k\in\mbb{N}$. We will prove it by induction on $k$.

When $k=0$, (3.66) holds by the definition (3.53). Consider $k=1$.
Write $\mu=\sum_{i=1}^n\mu_i\ves_i\in\Lmd^+$ with ${\cal
S}(\mu)=\{n_0,n_1,...,n_s\}$.  According to Lemma 3.2 and Lemma 3.4
with $M=V(\mu)$, the eigenvalues of $\vf|_{\widehat{V(\mu)}_{\la
1\ra}}$ are among
$$\{b-n,\;b+\mu_{n_{i-1}}-n_{i-1},\;b-\mu_{n_i}-2n+n_i\;\;\for\;\;i\in\ol{1,s}\}.\eqno(3.67)$$
 Recall that $\mu_r\in\mbb{N}/2$ for
$r\in\ol{1,n}$,
$$\mu_{\iota+1}-\mu_\iota\in\mbb{N}\;\;\for\;\;\iota\in\ol{1,n-1}\eqno(3.68)$$
and (2.12) holds. So
$$-\mu_{1+n_{i-1}}+n_{i-1},\;\mu_{n_i}+2n-n_i\in
\mu_1+2n-n_1-\mbb{N}\;\;\for\;\;i\in\ol{1,s}.\eqno(3.69)$$ If
$b\not\in \mu_1+2n-n_1-\mbb{N}$ and $b\neq n$, then all the
eigenvalues of $\vf|_{\widehat{V(\mu)}_{\la 1\ra}}$ are nonzero.  In
the case $\mu=(\sum_{i=1}^n\ves_i)/2$, the eigenvalues
$\vf|_{\widehat{V(\mu)}_{\la 1\ra}}$ are $b-n$ and $b+1/2$, which
are not equal to 0 because of $b\not\in n-\mbb{N}/2$. Thus (3.66)
holds for $k=1$.

Suppose that (3.66) holds for $k\leq \ell$ with $\ell\geq 1$.
Consider $k=\ell+1$. Note that
$$\vf(\widehat{V(\mu)}_{\la \ell+1\ra})=\sum_{i=0}^{2n}\vf(x_i\widehat{V(\mu)}_{\la \ell\ra})
=\sum_{i=0}^{2n}J_i[\vf(\widehat{V(\mu)}_{\la
\ell\ra})]=\sum_{i=0}^{2n}J_i(\widehat{V(\mu)}_{\la
\ell\ra})\eqno(3.70)$$ by the inductional assumption.  To prove
(3.66) with $k=\ell+1$ is equivalent to prove
$$\sum_{i=0}^{2n}J_i(\widehat{V(\mu)}_{\la
\ell\ra})=\widehat{V(\mu)}_{\la \ell+1\ra}.\eqno(3.71)$$

For any $u\in \widehat{V(\mu)}_{\la \ell-1\ra}$, Lemma 3.5 says that
$$J_0(x_0u)+\sum_{i=1}^n[J_i(x_{n+i}u)+J_{n+i}(x_iu)]=(2b-2n+\ell)\eta
u.\eqno(3.72)$$ Since $b\not\in n-\mbb{N}/2$, we have
$2b-2n+\ell\neq 0$, and so (3.72) gives
$$\eta u\in \sum_{i=0}^{2n}J_i(\widehat{V(\mu)}_{\la
\ell\ra})\qquad\for\;\;u\in  \widehat{V(\mu)}_{\la
\ell-1\ra}.\eqno(3.73)$$

Let $g\otimes v\in \widehat{V(\mu)}_{\la \ell\ra}$. According to
(3.40)-(3.46) and Lemma 3.4,
$$J_0(g\otimes v)=-\eta\ptl_{x_0}(g)\otimes v+g[(\ell+b+\td\omega)(x_0\otimes
v)],\eqno(3.74)$$
$$ J_i(g\otimes v)=-\eta\ptl_{x_{n+i}}(g)\otimes v
+g[(\ell+b+\td\omega)(x_i\otimes v)],\eqno(3.75)$$
$$J_{n+i}(g\otimes v)=-\eta\ptl_{x_i}(g)\otimes v
+g[(\ell+b+\td\omega)(x_{n+i}\otimes v)]\eqno(3.76)$$ for
$i\in\ol{1,n}$ (cf. (2.90) and (2.91)).
 Since
$$\eta\ptl_{x_i}(g)\otimes v\in
\sum_{r=0}^{2n}J_r(\widehat{V(\mu)}_{\la
\ell\ra})\qquad\for\;\;i\in\ol{0,2n}\eqno(3.77)$$ by (3.73),
Expressions (3.74)-(3.76) show
$$g[(\ell+b+\td\omega)(x_i\otimes
v)]\in \sum_{r=0}^{2n}J_r(\widehat{V(\mu)}_{\la
\ell\ra})\qquad\for\;\;i\in\ol{0,2n},\;g\in{\cal
A}_\ell.\eqno(3.78)$$

According to Lemma 3.2 and 3.4, the eigenvalue of
$(\ell+b+\td\omega)|_{\widehat{V(\mu)}_{\la 1\ra}}$ are among
$$\{b+\ell-n,b+\ell+\mu_{1+n_{i-1}}-n_{i-1},\;b+\ell-\mu_{n_i}-2n+n_i\mid i\in\ol{1,s}\}.\eqno(3.79)$$
Again $$-\ell-\mu_{1+n_{i-1}}+n_{i-1},\;-\ell-\mu_{n_i}+2n-n_i\in
\mu_1+2n-n_1-\mbb{N}\;\;\for\;\;i\in\ol{1,s}.\eqno(3.80)$$ If
$b\not\in\{n-\mbb{N}/2,\mu_1+2n-n_1-\mbb{N}\}$, then all the
eigenvalues of $(\ell+b+\td\omega)|_{\widehat{V(\mu)}_{\la 1\ra}}$
are nonzero. In the case $\mu=(\sum_{i=1}^n\ves_i)/2$, the
eigenvalues of $(\ell+b+\td\omega)|_{\widehat{V(\mu)}_{\la 1\ra}}$
are $b+\ell-n$ and $b+\ell+1/2$, which are not equal to 0 because of
$b\not\in n-\mbb{N}/2$. Hence
$$(\ell+b+\td\omega)(\widehat{V(\mu)}_{\la 1\ra})=\widehat{V(\mu)}_{\la
1\ra}.\eqno(3.81)$$ By (3.78) and (3.81),
$$g(\widehat{V(\mu)}_{\la 1\ra})\subset \sum_{r=0}^{2n}J_r(\widehat{V(\mu)}_{\la
\ell\ra})\qquad\for\;\;g\in{\cal A}_\ell,\eqno(3.82)$$ equivalently,
(3.66) holds for $k=\ell+1$. By induction, (3.66) holds for any
$k\in\mbb{N}.\qquad\Box$\psp

We remark that the $o(2n+3,\mbb{C})$-module $\widehat{V(\mu)}$  is
$o(2n+1,\mbb{C})$-finite, that is, $\widehat{V(\mu)}$ is of $({\cal
G},{\cal K})$-type with ${\cal G}=o(2n+3,\mbb{C})$ and ${\cal
K}=o(2n+1,\mbb{C})$. Up to this stage, we do not known if the
condition in Theorem 3.6 is necessary for the generalized conformal
$o(2n+3,\mbb{C})$-module $\widehat{V(\mu)}$ to be irreducible if
$\mu\neq 0$. In the case $\mu=0$, the situation becomes clear. \psp

 {\bf Theorem 3.7}. {\it The generalized conformal
$o(2n+3,\mbb{C})$-module $\widehat{V(0)}$ is irreducible if and only
if $b\not\in-\mbb{N}$. When $b=0$, $\widehat{V(0)}$ is isomorphic to
the natural conformal $o(2n+3,\mbb{C})$-module ${\cal A}$, on which
$A(f)=\vt(A)(f)$ for $A\in o(2n+3,\mbb{C})$ and $f\in{\cal A}$ (cf.
(2.45)-(2.48), (3.35) and (3.36)). The subspace $\mbb{C}$ forms a
trivial $o(2n+3,\mbb{C})$-submodule of the conformal module ${\cal
A}$ and the quotient space ${\cal A}/\mbb{C}$ forms an irreducible
$o(2n+3,\mbb{C})$-module.}

{\it Proof}.  Pick $0\neq v_0\in V(0)$. Then $\widehat{V(0)}={\cal
A}\otimes v_0$. We only list some facts different from the proof of
Theorem 2.7.
 Since $V(0)$ is the trivial
$o(2n+1,\mbb{C})$-module, $$J_0(f\otimes
v_0)=[x_0(D+b)-\eta\ptl_{x_0}](f)\otimes v_0.\eqno(3.83)$$
 Recall (3.22) and
(3.49). Set
$${\cal H}_k=\{f\in{\cal
A}_k\mid\Dlt(f)=0\}\qquad\for\;\;k\in\mbb{N}.\eqno(3.84)$$ Then
${\cal H}_k\otimes v_0$ is an irreducible
$o(2n+1,\mbb{C})$-submodule also with the highest-weight vector
$x_1^k\otimes v_0$. Moreover,
$$[\Dlt,\eta]=1+2n+D.\eqno(3.85)$$Furthermore,
$$J_0(x_1^{\ell-r-1}x_{n+1}^r\otimes v_0)=(b+\ell-1)x_1^{\ell-r-1}x_0x_{n+1}^r\otimes v_0
\eqno(3.86)$$ and $x_1^{\ell-2r-1}x_0\in{\cal H}_{\ell-2r}$. By the
similar arguments as those in the Proof of Theorem 2.7, we can prove
the conclusion in Theorem 3.7.$\qquad\Box$\psp

\begin{center}{\Large \bf Acknowledgement}\end{center}

 Part of this work was done during the first
author's visit to The University of Sydney, under the financial
support from Prof. Ruibin Zhang's ARC research grant. We thank Prof.
Zhang for his invitation, hospitality and helpful mathematical
discussion.

\end{document}